\documentclass[12pt]{article}

\newcommand{\vekk}[1]{} 

\usepackage{amsmath}
\usepackage{graphicx,psfrag,epsf}
\usepackage{enumerate}
\usepackage{natbib}
\usepackage{url} 

\vekk{
\addtolength{\oddsidemargin}{-.5in}%
\addtolength{\evensidemargin}{-.5in}%
\addtolength{\textwidth}{1in}%
\addtolength{\textheight}{1.3in}%
\addtolength{\topmargin}{-.8in}%
}

\RequirePackage[colorlinks,citecolor=blue,urlcolor=blue]{hyperref}

\newcommand{\pr}{\operatorname{\text{P}}} 
\newcommand{\E}{{\operatorname{\text{E}}}} 

\newcommand{\univar}{\operatorname{\text{U}}}
\newcommand{\Pareto}{\operatorname{\text{Pareto}}}

\newcommand{\thetaml}{\hat{\theta}_{\text{\tiny ML}}}
\newcommand{\thetamu}{\hat{\theta}_{\text{\tiny MU}}}

\newcommand{\thetarb}{\hat{\theta}_{\text{\tiny RB}}}

\newcommand{\thetalv}{\hat{\theta}_{\text{\tiny LV}}}

\newcommand{\thetabp}{\hat{\theta}_{\text{\tiny Bp}}}
\newcommand{\thetabGM}{\hat{\theta}_{\text{\tiny B2}}}
\newcommand{\thetaOPT}{\hat{\theta}_{\text{\tiny B3}}}
\newcommand{\thetabGT}{\hat{\theta}_{\text{\tiny B1}}}
\newcommand{\thetah}{\hat{\theta}}
\newcommand{\thetaSC}{\hat{\theta}_{\text{\tiny SC}}}

\newcommand{\st}{\operatornamewithlimits{{\mid}}}

\newcommand{\be}[1]{\begin{equation}\label{eq#1}} 
\newcommand{\ee}{\end{equation}}

\newcommand{\gttitle}{Fiducial Symmetry in Action}

\begin{document}


  \title{\vspace{-12ex}\bf \gttitle}
  \author{Gunnar Taraldsen \hspace{.2cm}\\
  Department of Mathematical Sciences\\
  Norwegian University of Science and Technology}
  \maketitle 
\begin{abstract}
  Symmetry is key in classical and modern physics.
  A striking example is conservation of energy as a consequence of
  time-shift invariance from Noether's theorem.
  Symmetry is likewise a key element in statistics,
  which, as also physics,
  provide models for real world phenomena.
  Sufficiency, conditionality, and invariance are examples of
  basic principles.
  \citet{GaliliMeilijson16ExampleImprovableRaoa} and
  \citet{Mandel20ScaledUniformModel} illustrate the first
  two principles very nicely by considering the scaled uniform model.
  We illustrate the third principle by providing further results which
  give optimal inference
  for the scaled uniform by symmetry considerations.
  The proofs are simplified by relying on fiducial arguments as initiated
  by \citet{Fisher30InverseProbability}.
\end{abstract}

\noindent%
{\it Keywords:}
Data generating equation;
Optimal equivariant estimate;
Scale family;
Conditionality principle;
Minimal sufficient;
Uniform distribution;

\vfill

\tableofcontents

\vfill

\newpage
\section{Introduction}
\label{sec:intro}

Let
\be{Fid}
y = \theta u
\ee
If $u = (u_1, \ldots, u_n)$ is a random sample from the
uniform law $\univar[1-k, 1+k]$ with known design parameter $k \in (0,1)$,
then $y = (y_1, \ldots, y_n)$ is a random sample from
the scaled uniform distribution
$\univar [\theta (1-k), \theta (1+k)]$
with scale parameter $\theta > 0$.
Estimation of $\theta$ in this normalized case
was investigated in some detail by
\citet{GaliliMeilijson16ExampleImprovableRaoa} and
\citet{Mandel20ScaledUniformModel}.
The remainder of this section recaps some of their results.
Additional results and discussion are
found in many classical texts  in theoretical statistics
since the scaled uniform, and its relatives, are used as a prototypical counterexamples
to results depending on a smooth likelihood:
Fisher Information, Crámer-Rao bound, efficiency of MLEs, $\ldots$

\vekk{
The choice of parametrisation gives $\E U_i = 1$ and
$\E Y_i = \E (\theta U_i)= \theta$,
so the scale parameter equals the mean.
Unbiased estimation of $\theta$ is in particular possible.
In the following main sections we will provide
novel and optimal inference
for the scaled uniform based on symmetry and fiducial arguments.
Additional results and discussion are
found in many classical texts  in theoretical statistics
since the scaled uniform is used as a prototypical counterexample
to results
depending on a smooth likelihood:
Fisher Information, Crámer-Rao bound, efficiency of MLEs, $\ldots$
}

The assumed data generating equation~(\ref{eqFid}) gives the likelihood

\be{Lik}
L (\theta) = \prod_i \frac{[\theta (1-k) \le y_i \le \theta (1+k)]}{2 k \theta}
= \frac{(\thetaml \le \theta \le \thetamu)}{(2 k \theta)^n}
\ee
where $\thetaml = y_{(n)}/(1+k)$ and $\thetamu = y_{(1)}/(1-k)$. 
The likelihood is hence zero for $\theta < \thetaml$ and for
$\theta > \thetamu$, and has jumps at $\thetaml$ and $\thetamu$.
The estimates $\thetaml$ and $\thetamu$ give deterministic
information since $\thetaml \le \theta \le \thetamu$
is always true.
Formally, $[\thetaml,\thetamu]$ is a $100\%$ confidence interval.
A minimal sufficient statistic 
is given by the smallest and largest observation
$(y_{(1)},y_{(n)})$,
by the sure interval $[\thetaml,\thetamu]$,
or equivalently by $s = (y_{(n)}, y_{(1)}/y_{(n)})$.

The maximum likelihood estimator (MLE) from equation~(\ref{eqLik}) is,
as the notation suggests, equal to the lower bound $\thetaml$.
It can be observed that if the closed interval $[1-k,1+k]$ is replaced
by the open interval $(1-k,1+k)$,
then this would give an example where the MLE does not exist.
The MLE is, in fact,
inefficient and biased as proved by \citet{GaliliMeilijson16ExampleImprovableRaoa}.
It is, furthermore, unreasonable since it totally ignores
the information provided by the smallest observation $y_{(1)}$.

An alternative unbiased estimator is
the Rao-Blackwellization of $Y_1$: 
\be{RB}
\thetarb = \E^\theta (Y_1 \st S=s) = \frac{y_{(1)}+y_{(n)} }{2} 
\ee
It is claimed, wrongly, by \citet[p.109]{GaliliMeilijson16ExampleImprovableRaoa}
that $(Y_1 \st S=s) \sim \univar [y_{(1)},y_{(n)}]$.
The claim is wrong since
the conditional law of $Y_1$ has point masses at both endpoints of the interval.
This exemplifies that distributions that are neither continuous nor discrete
appears naturally in probability calculations.
Equation~(\ref{eqRB}) is, nevertheless, correct.

The estimator $\thetarb$ is optimal in a related location model,
but not so in the given scale model.
\citet{GaliliMeilijson16ExampleImprovableRaoa}
prove that the Rao-Blackwell estimator is, in fact, uniformly improved by
$\thetalv = c^{-} y_{(1)} + c^{+} y_{(n)}$
\vekk{
\be{LV}
\thetalv = c^{-} y_{(1)} + c^{+} y_{(n)}
\ee
}
with $c^{\pm} = 0.5 (1 \pm k)/(1 + k^2 (n-1)/(n+1))$.
This ensures that $\thetalv$ is unbiased,
and has minimum variance in the class of linear
functions of  $y_{(1)}$ and  $y_{(n)}$.
The Rao-Blackwell estimator $\thetarb$
has, however, an advantage compared to the uniform improvement $\thetalv$:
It gives always a feasible estimate in the sense of belonging to the
sure interval $[\thetaml,\thetamu]$. 
The Rao-Blackwell estimator $\thetarb$ is 
optimal in the class of feasible linear unbiased estimators.

In the quest for better estimators
\citet{GaliliMeilijson16ExampleImprovableRaoa} 
turn to a Bayesian approach with an improper prior
$\pi (\theta) = \theta^{-p}$.
This is  a natural choice since combined with
the likelihood in equation~(\ref{eqLik}) the resulting posterior
is a truncated $\Pareto (\alpha, [a, b])$ distribution
with density:
\be{Pareto}
\pi (\theta \st y) = \frac{\alpha}{a^{-\alpha} - b^{-\alpha}}
\cdot (a \le \theta \le b) \cdot \theta^{-\alpha - 1}
\ee
The truncation parameters is $[a,b] = [\thetaml,\thetamu]$ and
the index is $\alpha = n + p - 1$.
Let $b_* = \thetamu/\thetaml$ and assume $\alpha \neq 1$.
A family of Bayes estimators is then
\be{Bp}
\thetabp = \E (\Theta \st Y=y)
= \frac{\alpha}{\alpha - 1} \cdot
\frac{1 - {b_*}^{1 - \alpha}}{1 - {b_*}^{-\alpha}}
\cdot \thetaml
\ee
%
\citet{GaliliMeilijson16ExampleImprovableRaoa} prove
that the estimator $\thetabGM$ is in fact unbiased.
Furthermore,
numerical evidence indicates that it is
uniformly better than $\thetalv$. 
It is quite remarkable that the Bayes estimator
from the prior $\pi (\theta) = \theta^{-2}$ is both unbiased and so good.

In enters \citet{Mandel20ScaledUniformModel} and gives an
alternative justification for the
estimator $\thetabGM$.
He bases the inference on the minimal sufficient
$s = (y_{(n)}, y_{(1)}/y_{(n)})$.
The key observation is that $s_2$ is ancillary:
The law of
$S_2 = Y_{(1)}/Y_{(n)} = U_{(1)}/U_{(n)}$ does not depend on $\theta$.
It is tempting, but erroneous, to conclude
from this that all information regarding $\theta$
is included in $S_1 = Y_{(n)}$.
As explained earlier, $\thetaml$, or equivalently $S_1$,
does not contain all information regarding $\theta$.

The conditionality principle dictates,
as advocated by \citet{Mandel20ScaledUniformModel},
that inference should be based on the
conditional model given the ancillary
$S_2 = s_2$.
Calculus shows that the
conditional law of $S_1$ given $S_2 = s_2$
is $\Pareto (\alpha, [a,b])$
with truncation interval
$[a,b] = \theta [(1-k)/s_2, (1+k)]$
and index $\alpha = -n$
\citep[eq.1]{Mandel20ScaledUniformModel}.
The conditional model has $\theta$ as a scale parameter.
A calculation using the expectation of the
$\Pareto (\alpha, [a,b])$ gives
$\E (S_1 \st S_2 = s_2) = \theta / \phi(s_2)$.
This gives an unbiased estimator $T = \phi(S_2) S_1$ for $\theta$ given $S_2 = s_2$.
It follows then that $T$ is also
unconditionally unbiased since
$\E (T) = \E (\E(T \st S_2)) = \E (\theta) = \theta$.
Remarkably, calculus shows that $T=\thetabGM$.

The arguments of \citet{Mandel20ScaledUniformModel}
give hence a simplified proof of unbiasedness of $\thetabGM$,
and identifies it with a natural estimator arising from
a conditional argument.
A natural question next is:
Is $\thetabGM$ optimal?
This is answered in the negative by
\citet{Mandel20ScaledUniformModel}.
He shows that $\thetabGM$ fulfills criteria
that ensures existence of an unbiased estimator
with smaller variance at any fixed $\theta_0$.
The next sections provide alternative estimators that are optimal.

\section{Fiducial Inference}
\label{sFid}

\citet[p.532,p.534]{Fisher30InverseProbability} introduced
fiducial inference for the correlation coefficient $\rho$
based on the empirical correlation $r$ of a random sample from
a two-dimensional Gaussian distribution. 
The argument is based on inversion of the
pivotal $F (R \st \rho) \sim \univar [0,1]$ where
$F (r \st \rho) = \pr (R \le r)$.
A draw from the fiducial distribution of the correlation coefficient
given $r$ is obtained by a draw $u \sim \univar [0,1]$, 
and then returning the unique solution $\rho$ of
the equation $F (r,\rho) = u$.

Unfortunately, Fisher gave many other recipes for obtaining
a fiducial which where conflicting, and the
birth of fiducial inference was a thorny one
\citep[p.185-]{SchwederHjort16ConfidenceLikelihoodProbabilitya}.
A particular blind alley is given by Fisher's
many attempts at using the likelihood $L$
instead of the cumulative $F$ for the construction of a fiducial.
Note that a particular likelihood can result from many different
statistical models and that a particular statistical model
can result from many different data generating equations.
\citet{CuiHannig19NonparametricGeneralizedFiducial}
give further references to the literature on fiducial inference,
and demonstrate that the initial approach of  
\citet{Fisher30InverseProbability} was the correct one
even in a non-parametric setting.

Modern fiducial inference is based on a data generating equation
instead of a family of distributions as a model for the observed data.
A given family of distributions can be the result of many
different data generating equations.
A particular data generating equation contains hence information
beyond the resulting family of distributions.
The data generating equation represents prior information,
but not in the form of a prior distribution as in Bayesian inference.
The fiducial argument was invented by Fisher to obtain
a distribution for the unknown parameter from
the observations and a model for cases where a Bayesian prior is absent.
It will next be demonstrated how this can be done for the
scaled uniform.
Furthermore, it will be proved in Section~\ref{sOpt} by equation~(\ref{eqFreqRisk})
that this fiducial gives an estimate with uniformly minimum
expected squared error.

A fiducial  model is specified by taking 
equation~(\ref{eqFid}), $y = \theta u$, as 
a data generating equation for the scaled uniform.
The direct recipe from the previous
is then to draw $u$ and solve for the parameter $\theta$.
Unfortunately, this fails since there are $n$ equations
for the one unknown $\theta$.
Intuitively, the solution is given by observing that
$u$ must be conditioned to be a point on the ray
$l_y = \{u \st u = \theta^{-1} y , \theta > 0\}$.
Mathematically, this is not sufficient, but it is
by instead conditioning on a function of $u$ with this line as a level set.
The set of different possible rays gives a partition of
the set of possible $u$'s.
Let $m(u)$ have exactly the partition as its level sets.
The required conditioning is then uniquely determined
by conditioning $u$ on $m (u) = m_y$ where
$l_y = \{u \st m(u) = m_y\}$.
The recipe is then to draw $u$ from its
initial distribution conditionally on $m (u) = m_y$,
and then return $\theta = y_n / u_n$ as a unique draw from the fiducial.

The previous argument determines the fiducial uniquely
from the data generating equation~(\ref{eqFid}),
and in fact for any scale model.
The analysis is simplified for the scaled uniform
by replacing equation~(\ref{eqFid}) by
the sufficient data generating equation
\be{SuffFid}
(y_{(1)}, y_{(n)}) = \theta (u_{(1)}, u_{(n)})
\ee
The fiducial is
$\Theta = y_{(n)} /V$
where
$V \sim (U_{(n)} \st U_{(1)}/U_{(n)} = s_2)$.
This is as in the previous recipe
since the mapping $(U_{(1)},U_{(n)})  \mapsto U_{(1)}/U_{(n)}$
has exactly the rays as its level sets. 
The calculations of \citet[eq.1]{Mandel20ScaledUniformModel} gives
now that $V \sim \Pareto (\alpha, [a,b])$
with truncation interval $[a,b] = [(1-k)/s_2, (1+k)]$
and index $\alpha = -n$.
The resulting fiducial is
\be{ThetaF}
\Theta \sim \Pareto (n, [\thetaml,\thetamu])
\ee
It has here been used that
$V \sim \Pareto (\alpha, [a,b])$ implies
$y_{(n)} V^{-1}$ $\sim$\\
$\Pareto (-\alpha, [y_{(n)}/b,y_{(n)}/a])$.

In a certain sense the analysis is now finished.
The fiducial in equation~(\ref{eqThetaF}) gives our state of knowledge
regarding the unknown model $\theta$ based on the fiducial model
in equation~(\ref{eqFid}),
or equivalently equation~(\ref{eqSuffFid}),
and the observation $y$.
This is, as Fisher intended, an alternative to a Bayesian posterior distribution.

\section{Optimal Decisions}
\label{sOpt}

The analysis can be continued in many ways depending
on the question of interest.
The original question here was:
How should the parameter $\theta$ be estimated?
This question is now simpler than initially because we have a probability
distribution for $\theta$ given in equation~(\ref{eqThetaF}).
What guess for $\theta$ should we choose when
its probability distribution is known?

Decision theory gives a possible route.
The guess, or estimate, is a decision.
One possibility is to minimize the risk given by
\be{risk1}
r = - \E \delta (\Theta - \thetah)
\ee
where $\delta$ is the Dirac delta function.
The density of $\Theta$ has a maximum at
$\thetaml$, so the result is then
the maximum likelihood estimator.
In a Baysian analysis this corresponds to choosing
the maximum a posterior (MAP) estimate.

Another possibility
is to minimize the risk given by
\be{risk2}
r = \E (\Theta - \thetah)^2
\ee
corresponding to expected squared error loss $(\theta - \thetah)^2$.
The resulting estimate is then
$\thetabGT$ as given in equation~(\ref{eqBp})
with $p=1$.
This is so because the fiducial coincides
with Bayesian posterior for the case $p=1$.
It follows similarly that 
the estimator $\thetabGM$ minimizes the risk
%
$r = \E [\Theta^{-1} (\Theta - \thetah)^2]$
%

Except for $\thetabGM$ none of the above estimators are unbiased.
The property of being unbiased is a natural demand
for estimating a location parameter,
but not so for a scale parameter.
A natural demand is to require scale invariance.
This translates into demanding that
$\thetah / \theta$ is a pivotal quantity.
An estimator $\thetah$ is said to be scale equivariant if this
demand is fulfilled.
As also noted by \citet[p.110]{GaliliMeilijson16ExampleImprovableRaoa}
this requirement holds for all suggested estimators here.
The natural question is then:
Is there an optimal scale equivariant estimator?
As explained by \citet[p.388-]{BERGER}
it is natural to only consider risk
corresponding to invariant loss.

A possible invariant risk is given by
\be{risk4}
r = \E (1 - \Theta^{-1} \thetah)^2
\ee
Differentiation and solving $0 = \partial_{\thetah} r$ gives
the optimal invariant estimator
\be{thetaOpt}
\thetaOPT = \frac{\E (\Theta^{-1})}{\E (\Theta^{-2})}
\ee
As the notation suggests this equals the Bayes posterior
estimator given in equation~(\ref{eqBp}) with $p=3$.
This can be seen by observing
$(1 - \thetah/\Theta)^2 = \theta^{-2} (\theta - \thetah)^2$,
and using the simple form of the fiducial distribution
in this particular case.

The observation also shows that $\thetaOPT = \psi (y)$
has uniformly  minimal squared error loss in the class
of scale equivariant estimators.
The claim follows from the frequentist risk
\be{FreqRisk}
\begin{split}
r^\theta & = \E^\theta (\theta - \psi (Y))^2 =
\theta^2 \E (1 - \theta^{-1} \psi(\theta U))^2 \\
& =
\theta^2 \E \left[ \E \left\{(1 - \Theta^{-1 } \psi(\Theta U))^2 \st S_2 \right\}\right] =
\theta^2 \E \left[ \E \left\{(1 - \Theta^{-1 } \psi(y))^2 \st S_2 \right\}\right]
\end{split}
\ee
\citet{TaraldsenLindqvist13fidopt} prove 
similar frequentist optimality more generally
for equivariant decision rules derived from
an invariant loss and a fiducial distribution. 

Another possible invariant risk is given by
\be{risk5}
r = \E (\ln \Theta - \ln \thetah)^2
\ee
The corresponding optimal equivariant estimator
$\thetaSC$ is determined by $\ln \thetaSC = \E \ln \Theta$.
Equation~(\ref{eqThetaF}) gives
$\Theta \sim \thetaml \Pareto(n, [1,b_*])$
with $b_* = \thetamu/\thetaml$ and calculus gives
\be{OptLog}
\thetaSC = \thetaml
\exp \left\{
  \frac{1}{n} -
  \frac{\ln b_*}{1-{b_*}^n}
\right\}
\ee

\section{Conclusion}
\label{sConc}

It has been proven that the estimators
$\thetaOPT$ and $\thetaSC$ given in equation~(\ref{eqBp})
and equation~(\ref{eqOptLog}) respectively
are optimal equivariant estimators
in the sense of uniformly
minimizing the frequentist
risks
$\E^\theta (\thetaOPT - \theta)^2$ and
$\E^\theta (\ln \thetaSC - \ln \theta)^2$
respectively.
These results follows from invariance and
complements the results for
the scaled uniform as previously discussed by 
\citet{GaliliMeilijson16ExampleImprovableRaoa} and
\citet{Mandel20ScaledUniformModel}.
Sufficiency, conditionality, and invariance are central
themes in theoretical statistics and in the arguments given here.

Existence of a uniformly minimum-variance
unbiased (UMVU) estimator is left open.
It is known,
as demonstrated by \citet[p.379, Exercise 11]{RAO},
that in the case where the minimal sufficient statistic $S$ is not complete
some parameters may have UMVU estimators,
but others may not.
It would be interesting,
based on curiosity,
to know if a UMVU estimator for $\theta$ exists,
but equivariance is a more natural demand for the case of a scale parameter.

In retrospect, it was perhaps fortunate that Fisher invented many and conflicting
roads in his attempts to arrive at fiducial inference.
Blind alleys are blind alleys, but often there are rewards along the way.
This has certainly been the case also for the many different versions
of fiducial inference that has been suggested.
\citet[p.532,p.534]{Fisher30InverseProbability} invented both confidence intervals
and confidence distributions in his initial attempt.
The theory of the latter topic is still in its infancy,
but considerable progress have been made lately as documented by 
\citet{SchwederHjort16ConfidenceLikelihoodProbabilitya}.

The fiducial distribution of $\Theta$ is
$\Pareto (n, [\thetaml, \thetamu])$ and it is a confidence distribution.
This follows from the general arguments given by 
\citet{TaraldsenLindqvist13fidopt}. 
It can be concluded that fiducial inference as
initiated by \citet{Fisher30InverseProbability}
can be used to obtain both optimal estimators and
exact confidence distributions.
This has here been demonstrated for the
scaled uniform distribution
which is otherwise used as a counterexample
for obtaining good theoretical results.

\bibliographystyle{Chicago}

\bibliography{tex}
\end{document}